\documentclass[onecolumn,12pt]{autart}
\usepackage{graphicx}
\input{amssym}

\newcommand{\RR}{{\Bbb R}}
\renewcommand{\theequation}{\arabic{section}.\arabic{equation}}
\begin{document}
\begin{frontmatter}
\title{Projective analysis and preliminary group classification\\
of the nonlinear fin equation $u_t=(E(u)u_x)_x + h(x)u$}
\thanks[footnoteinfo]{ Corresponding author.}
\author[]{M. Nadjafikhah\thanksref{footnoteinfo}}\ead{m\_nadjafikhah@ius.ac.ir},
\author[]{A. Mahdipour--Shirayeh}\ead{mahdipour@iust.ac.ir},
\address{School of Mathematics, Iran University of Science and Technology,\\ Narmak, Tehran 1684613114, Iran.}
\renewcommand{\sectionmark}[1]{}
\begin{abstract}
In this paper we investigate for further symmetry properties of
the nonlinear fin equations of the general form $u_t=(E(u)u_x)_x +
h(x)u$ rather than recent works on these equations. At first, we
study the projective (fiber--preserving) symmetry to show that
equations of the above class can not be reduced to linear
equations. Then we determine an equivalence classification which
admits an extension by one dimension of the principal Lie algebra
of the equation. The invariant solutions of equivalence
transformations and classification of nonlinear fin equations
among with additional operators are also given.
\end{abstract}
\end{frontmatter}
\section{Introduction}
Investigations for symmetry properties of mathematical models of
heat conductivity and diffusion processes \cite{Ib1} are
traditionally formulated in terms of nonlinear differential
equations which often envisage us with difficulties in studying.
To solve this problem, symmetry methods play a key role for
finding their exact solutions, similar solutions
\cite{BK,Do,Pa,PS,VP} and invariants.
\\[3mm] \mbox{ }\hspace{3mm} In this study we generalize
the study of a class of nonlinear fin equations which has been
recently studied in some references and specially in \cite{PS,VP}.
So, we are dealing with the class of nonlinear fin equations of
the general form
\begin{eqnarray}
u_t=(E(u)\,u_x)_x + h(x)\,u, \label{eq:1}
\end{eqnarray}
in which we assumed that $E_u \neq 0$, $u$ is treated as the
dimensionless temperature, $t$ and $x$ the dimensionless time and
space variables, $E$ the thermal conductivity, $h = -N^2\,f(x)$,
$N$ the fin parameter and $f$ the heat transfer coefficient
\cite{BK}.
\\[3mm] \mbox{ }\hspace{3mm} The Lie point symmetry in linear and
nonlinear case, the condition $E_u=0$ corresponds to the linear
case, the class of nonlinear one-dimensional diffusion equations
when $h=0$, the class of diffusion–-reaction equations when
$h=\verb"cons."$, the case which the thermal conductivity is a
power function of the temperature and additional equivalence
transformations, conditional equivalence groups and nonclassical
symmetries have all investigated and listed in Table 1 of
\cite{VP}. The point symmetry group of nonlinear fin equations of
class (\ref{eq:1}) were considered in a number of papers (e.g. see
\cite{BK} for the physical meaning and applications of the
equation). The Lie algebra of the point symmetry group of
Eq.~(\ref{eq:1}) is
\begin{eqnarray}
{\goth g}_1:=\Big\langle\,\frac{\partial}{\partial
t}\,\Big\rangle, \label{eq:1-1}
\end{eqnarray}
In the next section, we concern with the problem of finding
projective symmetry group of Eq.~(\ref{eq:1}) as a special case of
the point symmetry group; since it may have important information
about the equation.
\\[3mm] \mbox{ }\hspace{3mm}   The equivalence classification of
Eq.~(\ref{eq:1}) in the special case
\begin{eqnarray}
u_t=(E(u)\,u_x)_x,
\end{eqnarray}
has performed by L.V.~Ovsiannikov \cite{Ov}. In \cite{PS} authors
carried out the another special case of equivalence group
\begin{eqnarray}
\widetilde{t}=\delta_1\,t+\delta_2,\hspace{1cm} \widetilde{x} =
\delta_3\,x + \delta_4, \hspace{1cm} \widetilde{u} = \delta_5\,u,
\hspace{1cm} \widetilde{E} = \delta_1^{-1} \delta_3^2\,E,
\hspace{1cm} \widetilde{h} = \delta^{-1}_1\,h,
\end{eqnarray}
of Eq.~(\ref{eq:1}) when $\delta_i,\; i=1, \cdots, 5$, are
arbitrary constants and $\delta_1\,\delta_3\,\delta_5 \neq 0$. The
more general class of nonlinear fin equations is the nonlinear
heat conductivity equations of the form
\begin{eqnarray}
u_t = F(t, x, u,u_x)u_{xx} + G(t, x, u, u_x ), \label{eq:1-2}
\end{eqnarray}
which admits non-trivial symmetry group. The group classification
of (\ref{eq:1-2}) is presented in some references \cite{BLZ,LS}.

However, since the equivalence group of (\ref{eq:1-2}) is
essentially wider than those for particular cases, the results of
\cite{BLZ,LS} cannot be directly used for symmetry classification
of particular ones. Nevertheless, these results are useful for
finding additional equivalence transformations in the class of our
problem. Therefore in contrast to the above works, in the last two
sections of this paper, we study group classification of
Eq.~(\ref{eq:1}) under equivalence transformations in the general
case. Furthermore, a number of nonlinear invariant models which
have nontrivial invariance algebras are obtain.
\\[3mm] \mbox{ }\hspace{3mm}    From \cite{ZL} we know that if the
partial differential equation possesses non-trivial symmetry, then
it is invariant under some finite-dimensional Lie algebra of
differential operators which is completely determined by its
structural constants. In the event that the maximal algebra of
invariance is infinite--dimensional, then it contains, as a rule,
some finite-dimensional Lie algebra. Also, if there are local
non-singular changes of variables which transform a given
differential equation into another, then the finite-dimensional
Lie algebra of invariance of these equations are isomorphic, and
in the group-theoretic analysis of differential equations such
equations are considered to be equivalent. To realize the group
classification, we use of the proposed approach consists in the
implementation of an algorithm explained and performed in
references \cite{BLZ,Ib3,Ov,PNI}. For this goal, our method is
completely similar to the way of \cite{LH} for the nonlinear wave
equation $u_{tt} = f(x,u)u_{xx} + g(x,u)$.
\section{Projective symmetries of nonlinear fin equations}
In this section, we are concerning with group classification of
nonlinear fin equations by projective transformations group as a
special case of the point symmetry group. Our study is based on
the method of \cite{Ol} for Lie infinitesimal method.
\\[3mm] \mbox{ }\hspace{3mm}    The equation is a relation among
with the variables of 2--jet space $J^2(\RR^2,\RR)$ with (local)
coordinate
\begin{eqnarray}
(t,x,u,u_t,u_x,u_{tt},u_{tx},u_{xx}),\label{var}
\end{eqnarray}
 where this coordinate involving independent variables $t,x$ and dependent variable $u$
and derivatives of $u$ in respect to $t$ and $x$ up to order 2
(each index will indicate the derivation with respect to it,
unless we specially state otherwise). Let ${\mathcal M}$ be the
total space of independent and dependent variables resp. $t,x$ and
$u$. The solution space of Eq.~(\ref{eq:1}), (if it exists) is a
subvariety $S_{\Delta}\subset J^2(\RR^2,\RR)$ of the second order
jet bundle of 2-dimensional sub-manifolds of ${\mathcal M}$. If we
wish to preserve the bundle structure of the space ${\mathcal M}$,
we must restrict to the class of fiber--preserving transformations
in which the changes in the independent variable are unaffected by
the dependent variable. {\it Projective} or {\it fiber--preserving
symmetry group} on ${\mathcal M}$ is introduced by transformations
in the form of
\begin{eqnarray}
\tilde{x}=\phi(t,x),\hspace{1cm}\tilde{t}=\chi(t,x),\hspace{1cm}
\tilde{u}=\psi(t,x,u),\label{eq:2}
\end{eqnarray}
for arbitrary smooth functions $\phi,\chi,\psi$. Also assume that
\begin{eqnarray}
v:=\xi^1(t,x)\,\displaystyle{\frac{\partial }{\partial
t}}+\xi^2(t,x)\,\displaystyle{\frac{\partial }{\partial x}}+
\eta(t,x,u)\,\displaystyle{\frac{\partial }{\partial u}},
\label{eq:2-1}
\end{eqnarray}
with coefficients as arbitrary smooth functions be the general
form of infinitesimal generators which signify the Lie algebra
${\frak g}$ of the projective symmetry group $G$ of
Eq.~(\ref{eq:1}). The second order prolongation of $v$
\cite{Ol,Ov} as a vector field on $J^2(\RR^2,\RR)$ is as follows
\begin{eqnarray}
v^{(2)}:=v + \eta^t\,\displaystyle{\frac{\partial }{\partial
u_t}}+\eta^x\,\displaystyle{\frac{\partial }{\partial u_x}}+
\eta^{tt}\,\displaystyle{\frac{\partial }{\partial u_{tt}}}+
\eta^{tx}\,\displaystyle{\frac{\partial }{\partial u_{tx}}}+
\eta^{xx}\,\displaystyle{\frac{\partial }{\partial u_{xx}}},
\end{eqnarray}
where $\eta^t,\eta^x$ and $\eta^{tt},\eta^{tx},\eta^{xx}$ are
arbitrary smooth functions depend to variables $t,x,u,p,q$ and
(\ref{var}) resp. These coefficients are computed as following
\begin{eqnarray}
\eta^J &=& {\mathcal D}_J(Q) + \xi^1\,u_{J,t}+ \xi^2\,u_{J,x}
\end{eqnarray}
where ${\mathcal D}$ is total derivative, $J$ is a multi-index
with length $1\leq |J|\leq2$ of variables $t,x$ and $Q= u -
\xi^1\,u_t - \xi^2\,u_x$ is characteristic of $v$ \cite{Ol}.
According to \cite{Ol}, $v$ is a projective infinitesimal
generator of Eq.~(\ref{eq:1}) if and only if
$v^{(2)}[\mbox{Eq.~(\ref{eq:1})}]=0$. By applying $v^{(2)}$ on the
equation we have the following equation
\begin{eqnarray}
\xi^2\,h_x\,u+\eta(D_{uu}\,q^2+D_u\,q_x)-\eta^t +
2\,\eta^x\,D_u\,q+\eta^{xx}\,D=0,\hspace{1cm}\mbox{whenever}\hspace{0.5cm}
\mbox{Eq.~(\ref{eq:1}) is satisfied}.
\end{eqnarray}
In the extended form of the latter equation, functions
$\xi^1,\xi^2$ and $\eta$ only depend to $t, x, u$ rather than
other variables, i.e. $u_t, u_x, u_{tt}, u_{tx}, u_{xx}$, hence
the equation will be satisfied if and only if the individual
coefficients of the powers of $u_t,u_x$ and their multiplications
vanish. This tends to the following over-determined system of
determining equations
\begin{eqnarray}
&&E_u\,\xi^1_x=0, \hspace{1cm}   E\,\xi^1_x=0, \hspace{1cm} \eta_u-\xi^1_t+E\,\xi^1_{xx}=0,  \nonumber\\
&&E_u\,\eta + E\,\eta_u-2\,E\,\xi^2_x=0, \hspace{1cm}
\xi^2\,h_x\,u-\eta_t+E\,\eta_{xx}+h\,\eta=0, \label{eq:4}\\
&&E_{uu}\,\eta + 2\,E_u(\eta_u-\xi^2_x)+E\,\eta_{uu}=0,
\hspace{1cm} \xi^2_t+2\,E_u\,\eta_x+E(2\,\eta_{xu}-\xi^2_{xx})=0.
\nonumber
\end{eqnarray}
By solving Eq.~(\ref{eq:4}), the general solution to these
differential equations for $\xi^1,\xi^2$ and $\eta$ will be found:
\begin{eqnarray}
\xi^1(t,x)=c, \hspace{1cm} \xi^2(t,x)=0, \hspace{1cm}
\eta(t,x,u)=0,
\end{eqnarray}
with arbitrary constant $c$. Therefore the Lie algebra ${\frak g}$
spanned by projective infinitesimal generators of Eq.~(\ref{eq:1})
is ${\frak g}=\langle \frac{\partial}{\partial t}\rangle$ and the
projective symmetry group is nothing but the time translation
group.
{\thm A complete set of all infinitesimal generators of nonlinear
fin equation up to projective transformations admits the structure
of one-dimensional Lie algebra and so is isomorphic to $\RR$. The
point symmetry group and projective symmetry of nonlinear fin
equation are equal.}
\\[3mm] \mbox{ }\hspace{3mm}    It is well--known that the existence
of a non--fiber--preserving symmetry usually indicates that one
can significantly simplify the equation by some kind of
hodograph--like transformation interchanging the independent and
dependent variables. But in the case of our problem, we can not
use this advantage for simplifying nonlinear fin equations.
\\[3mm] \mbox{ }\hspace{3mm}    According to \cite{Ol}, any system
of partial differential equations which has only a
finite--dimensional symmetry group is certainly not linearizable,
that is, for every change of variables, it can not be mapped to an
inhomogeneous form of the linear system ${\mathcal D}[u]=f$, where
${\mathcal D}$ is a second order linear differential operator, $u$
indicates dependent variables and $f$ denotes smooth functions of
independent variables.

{\rem A system of nonlinear fin equations in the form of
Eq.~(\ref{eq:1}) can not be reduced into an inhomogeneous form of
a linear system.}
\section{Equivalence transformations}
At the present section we follow the method of L.V.~Ovsiannikov
\cite{Ov} for partial differential equations. His approach is
based on the concept of an equivalence group, which is a Lie
transformation group acting in the extended space of independent
variables, functions and their derivatives, and preserving the
class of partial differential equations under study. It is
possible to modify Lie's algorithm in order to make it applicable
for the computation of this group \cite{BLZ,Ib3,Ov,PNI}. Next we
construct the optimal system of subgroups of the equivalence
group.
\\[3mm] \mbox{ }\hspace{3mm}    An {\it equivalence transformation}
is a non-degenerate change of the variables $t, x, u$ taking any
equation of the form (\ref{eq:1}) into an equation of the same
form, generally speaking, with different $E(u)$ and $h(x)$. The
set of all equivalence transformations forms an equivalence group
$G$. We shall find a continuous subgroup $G_C$ of it making use of
the infinitesimal method.
\\[3mm] \mbox{ }\hspace{3mm}    We investigate for an operator of
the group $G_C$ in the general form
\begin{eqnarray}
Y:=\xi^1(t,x)\,\displaystyle{\frac{\partial }{\partial
t}}+\xi^2(t,x)\,\displaystyle{\frac{\partial }{\partial x}}+
\eta(t,x,u)\,\displaystyle{\frac{\partial }{\partial u}} +
\varphi(t,x,u,E,h)\,\displaystyle{\frac{\partial }{\partial E}} +
\chi(t,x,u,E,h)\,\displaystyle{\frac{\partial }{\partial h}}.
\label{eq:10}
\end{eqnarray}
from the invariance conditions of Eq.~(\ref{eq:1}) written as the
system
\begin{eqnarray}
u_t = (E(u)\,u_x)_x + h(x)\,u, \hspace{1cm}  E_t = E_x =0,
\hspace{1cm} h_t = h_u = 0,\label{eq:11}
\end{eqnarray}
where we assumed that $u, E, h$ are differential variables: $u$ on
the base space $(t, x)$ and $E , h$ on the total space $(t, x,
u)$. Also in Eq.~(\ref{eq:10}) the coefficients depend to
variables $t, x, u$ and the two last ones, in addition, depend to
$E, h$. The invariance conditions of the system (\ref{eq:11}) are
\begin{eqnarray}
\widetilde{Y}\,[u_t - (E(u)\,u_x)_x - h(x)\,u]=0, \hspace{1cm}
\widetilde{Y}\,[E_t]=\widetilde{Y}\,[E_x]=0, \hspace{1cm}
\widetilde{Y}\,[h_t] =\widetilde{Y}\,[h_u] = 0, \label{eq:12}
\end{eqnarray}
where
\begin{eqnarray}
\widetilde{Y}:= Y +  \eta^t\,\displaystyle{\frac{\partial
}{\partial u_t}}+\eta^x\,\displaystyle{\frac{\partial }{\partial
u_x}} + \eta^{tt}\,\displaystyle{\frac{\partial }{\partial
u_{tt}}}+ \eta^{tx}\,\displaystyle{\frac{\partial }{\partial
u_{tx}}}+ \eta^{xx}\,\displaystyle{\frac{\partial }{\partial
u_{xx}}}+ \varphi^t\,\displaystyle{\frac{\partial }{\partial E_t}}
+ \varphi^x\,\displaystyle{\frac{\partial }{\partial E_x}} +
\chi^t\,\displaystyle{\frac{\partial }{\partial h_t}} +
\chi^u\,\displaystyle{\frac{\partial }{\partial
h_u}}.\label{eq:12-1}
\end{eqnarray}
is the prolongation of the operator (\ref{eq:10}). Coefficients
$\eta^J$ for multi--index $J$ (with length $1\leq |J|\leq 2$) have
given in section 2 and by applying the prolongation procedure to
differential variables $E, h$ with independent variables $(t, x,
u)$ we have
\begin{eqnarray}
&& \varphi^t = \widetilde{{\mathcal D}}_t(\varphi) -
E_t\,\widetilde{{\mathcal D}}_t(\xi^1)- E_x\,\widetilde{{\mathcal
D}}_t(\xi^2) - E_u\,\widetilde{{\mathcal D}}_t(\eta)=
\widetilde{{\mathcal D}}_t(\varphi) - E_u\,\widetilde{{\mathcal D}}_t(\eta),  \nonumber\\
&& \varphi^x = \widetilde{{\mathcal D}}_x(\varphi) -
E_t\,\widetilde{{\mathcal D}}_x(\xi^1)- E_x\,\widetilde{{\mathcal
D}}_x(\xi^2) - E_u\,\widetilde{{\mathcal D}}_x(\eta) =
\widetilde{{\mathcal D}}_x(\varphi) - E_u\,\widetilde{{\mathcal D}}_x(\eta),
\nonumber\\[-2.5mm]
&& \\[-2mm]
&& \chi^t = \widetilde{{\mathcal D}}_t(\chi) -
h_t\,\widetilde{{\mathcal D}}_t(\xi^1)- h_x\,\widetilde{{\mathcal
D}}_t(\xi^2) - h_u\,\widetilde{{\mathcal D}}_t(\eta) =
\widetilde{{\mathcal D}}_t(\chi) - h_x\,\widetilde{{\mathcal
D}}_t(\xi^2),  \nonumber\\
&& \chi^u = \widetilde{{\mathcal D}}_u(\chi) -
h_t\,\widetilde{{\mathcal D}}_u(\xi^1)- h_x\,\widetilde{{\mathcal
D}}_u(\xi^2) - h_u\,\widetilde{{\mathcal D}}_u(\eta)
=\widetilde{{\mathcal D}}_u(\chi) - h_x\,\widetilde{{\mathcal
D}}_u(\xi^2), \nonumber
\end{eqnarray}
where in view of Eq.~(\ref{eq:11}) we have
\begin{eqnarray}
&& \widetilde{{\mathcal D}}_t := \displaystyle{\frac{\partial
}{\partial t}} + E_t\,\displaystyle{\frac{\partial }{\partial E}}
+ h_t\,\displaystyle{\frac{\partial}{\partial h}} =
\displaystyle{\frac{\partial
}{\partial t}}, \nonumber \\
&& \widetilde{{\mathcal D}}_x := \displaystyle{\frac{\partial
}{\partial x}} + E_x\,\displaystyle{\frac{\partial }{\partial E}}
+ h_x\,\displaystyle{\frac{\partial}{\partial
h}}=\displaystyle{\frac{\partial }{\partial x}} + h_x\,\displaystyle{\frac{\partial}{\partial h}}, \\
&& \widetilde{{\mathcal D}}_u := \displaystyle{\frac{\partial
}{\partial u}} + E_u\,\displaystyle{\frac{\partial }{\partial E}}
+ h_u\,\displaystyle{\frac{\partial}{\partial h}} =
\displaystyle{\frac{\partial }{\partial u}} +
E_u\,\displaystyle{\frac{\partial }{\partial E}}.\nonumber
\end{eqnarray}
Substituting (\ref{eq:12-1}) in (\ref{eq:12}) we tend to the
following system
\begin{eqnarray}
&& \varphi\,u_{xx} + h\,\eta + \chi\,u - \eta^t +
2\,E_u\,u_x\,\eta^x + \varphi^u\,u_x^2
+ E\,\eta^{xx}=0,  \label{eq:13}\\
&& \varphi^t=0,  \hspace{1cm}\varphi^x=0, \label{eq:13-1}\\
&& \chi^t=0,  \hspace{1cm}\chi^u=0. \label{eq:14}
\end{eqnarray}

Replacing relations $\eta^J$ (for multi--index $J$ with length
$1\leq|J|\leq 2$) and $\chi^t, \chi^x$ in
Eqs.~(\ref{eq:13})--(\ref{eq:14}) and then introducing the
relation $u_t= (E(u)\,u_x)_x + h u$ to eliminate $u_t$, we have
five relations which are called determining equations. The four
last ones are the determining equations associated with
Eqs.~(\ref{eq:13-1}), (\ref{eq:14}), i.e.,
\begin{eqnarray}
&& \varphi_t - E_u\,\eta_t=0, \hspace{1cm} \varphi_x +
h_x\,\varphi_h - E_u\,\eta_x=0, \hspace{1cm}  \label{eq:15} \chi_t
- h_x\,\xi^2_t=0,\hspace{1cm} \chi_u + E_u\,\chi_E -
h_x\,\xi^2_u=0.
\end{eqnarray}
But these relations must hold for every $E$ and $h$ and this fact
results in the following conditions
\begin{eqnarray}
\xi^2_t = \xi^2_u =0,  \hspace{1cm}  \eta_t= \eta_x=0,
\hspace{1cm} \varphi_t= \varphi_x = \varphi_h =0, \hspace{1cm}
\chi_t = \chi_u = \chi_E =0,
\end{eqnarray}
so, we find that
\begin{eqnarray}
\xi^2=\xi^2(x),  \hspace{1cm}  \eta=\eta(u), \hspace{1cm}
\varphi=\varphi(u, E), \hspace{1cm} \chi=\chi(x, h).
\end{eqnarray}
Adding these conditions to the remained determining equation,
since $u_t, u_x, u_{tt}, u_{tx}, u_{xx}$ are considered to be
independent variables, we lead to the following system of
equations
\begin{eqnarray}
&&h\,\eta + \chi\,u + h\,u\,(\xi^1_t- \eta_u )=0, \hspace{1cm} E_u\,(\eta_u + \xi^1_t -2\,\xi^2_x) + E\,\eta_{uu} + E_{uu}\,\eta =0, \nonumber\\
&&\\[-3mm] && \varphi + \varphi_u + E\,(\xi^1_t - 2\,\xi^2_x) =0,
\hspace{1cm} E\,\xi^1_x =0, \hspace{1cm} E\,\xi^1_u = 0,
\hspace{1cm} E\,\xi^2_{xx}=0.\nonumber
\end{eqnarray}
This system follows
\begin{eqnarray}
\xi^1 = 2\,c_1\,t + c_2, \hspace{1cm}\xi^2 = c_1\,x + c_3,
\hspace{1cm} \eta=0, \hspace{1cm} \varphi =e^{-u}\,F(E),
\hspace{1.1cm} \chi = -2\,c_1\,h,
\end{eqnarray}
with arbitrary function $F=F(E)$ and constants $c_1, c_2, c_3$.
Therefore the class of Eqs.~(\ref{eq:1}) has an infinite
continuous group of equivalence transformations generated by
infinitesimal operators
\begin{eqnarray}
&& Y_1=\displaystyle{\frac{\partial }{\partial t}}, \hspace{1cm}
Y_2=\displaystyle{\frac{\partial }{\partial x}}, \hspace{1cm}
Y_3=2\,t\,\displaystyle{\frac{\partial }{\partial t}}+
x\,\displaystyle{\frac{\partial }{\partial x}}-
2\,h\,\displaystyle{\frac{\partial }{\partial h}}, \hspace{1cm}
Y_F=e^{-u}\,F(E)\,\displaystyle{\frac{\partial }{\partial
E}}.\label{eq:15-1}
\end{eqnarray}
Moreover, in the group of equivalence transformations are included
also discrete transformations, i.e., reflections
\begin{eqnarray}
t\longmapsto -t,  \hspace{1cm} x \longmapsto -x, \hspace{1cm} u
\longmapsto -u, \hspace{1cm} E\longmapsto -E, \hspace{1cm}
h\longmapsto -h.
\end{eqnarray}

The communication relations between these vector fields is given
in Table 1. The Lie algebra ${\goth g}:=\langle\, Y_F, Y_i:
i=1,\cdots,3\,\rangle$ is solvable since the descending sequence
of derived subalgebras of ${\goth g}$: ${\goth g}\supset{\goth
g}^{(1)}=\langle\, 2\,Y_1, Y_2\,\rangle\supset{\goth
g}^{(2)}=\{0\}$, terminates with a null ideal. But its Killing
form: $K(v,w)={\rm tr}({\rm ad}(v)\circ{\rm ad}(w))=5\,a_3\,b_3$
for each $v=\sum_i v_i Y_i + v_4 Y_F$ and $w=\sum_j w_j Y_j + w_4
Y_F$ in ${\goth g}$ is degenerate and hence ${\goth g}$ is neither
semisimple nor simple.
{\thm Let $G_i$ be the one--parameter groups generated by the
$Y_i$, then we have
\begin{eqnarray}
\begin{array}{lll}
G_1: (t,x,u,E,h)\longmapsto (t+s,x,u,E,h), &&
G_2: (t,x,u,E,h)\longmapsto (t,x+s,u,E,h), \\
G_3: (t,x,u,E,h)\longmapsto (t\,e^{2s}, x\,e^s, u, E, h\,e^{2s}),
&& G_4: (t,x,u,E,h)\longmapsto (t,x,u,\bar{E},h),
\end{array}
\end{eqnarray}
when $\bar{E}$ is the solution of equation
$\int^{\bar{E}}_c\,d\alpha /F(\alpha)=e^{-u}\,s$ for a constant
$c$. Furthermore, if $u=f(t,x)$ for functions $E$ and $h$ be a
solution of nonlinear fin equation, so are the following functions
\begin{eqnarray}
&& u_1=f(t+s,x,u),\hspace{1cm} u_2=f(t,x+s,u),
\end{eqnarray}
for the same functions $E$ and $h$, $u_3=f(t\,e^{2s}, x\,e^s, u)$
for the same functions $E$ and $\bar{h}=h\,e^{2s}$ and also
$u_4=f(t,x,u)$ for $\bar{E}$ and the same $h$.}
\section{Preliminary group classification}
In many applications of group analysis, most of extensions of the
principal Lie algebra admitted by the equation under consideration
are taken from the equivalence algebra ${\goth g}_{{\mathcal E}}$.
These extensions are called ${\mathcal E}$--extensions of the
principal Lie algebra. The classification of all nonequivalent
equations (with respect to a given equivalence group $G_{\mathcal
E}$) admitting ${\mathcal E}$--extensions of the principal Lie
algebra is called a {\it preliminary group classification}
\cite{Ib3}. We can take any finite--dimensional subalgebra
(desirable as large as possible) of an infinite--dimensional
algebra with basis (\ref{eq:15-1}) and use it for a preliminary
group classification. We select the subalgebra ${\goth g}_4$
spanned on the following operators:
\begin{eqnarray}
&& Y_1=\displaystyle{\frac{\partial }{\partial t}}, \hspace{1cm}
Y_2=\displaystyle{\frac{\partial }{\partial x}}, \hspace{1cm}
Y_3=2\,t\,\displaystyle{\frac{\partial }{\partial t}}+
x\,\displaystyle{\frac{\partial }{\partial x}}-
2\,h\,\displaystyle{\frac{\partial }{\partial h}}, \hspace{1cm}
Y_4=e^{-u}\,E\,\displaystyle{\frac{\partial }{\partial E}}.
\label{eq:15-2}
\end{eqnarray}

It is well-known that the problem of classifying invariant
solutions is equivalent to the problem of classifying subgroups of
the full symmetry group under conjugation in which itself is
equivalent to determining all conjugate subalgebras \cite{Ol,Ov}.
The latter problem, tends to determine a list (that is called an
{\it optimal system}) of conjugacy inequivalent subalgebras with
the property that any other subalgebra is equivalent to a unique
member of the list under some element of the adjoint
representation i.e. $\bar{{\goth h}}\,{\rm Ad}(g)\,{\goth h}$ for
some $g$ of a considered Lie group. Thus we will deal with the
construction of the optimal system of subalgebras of ${\goth
g}_4$.
\\[3mm] \mbox{ }\hspace{3mm}    The adjoint action is given by the
Lie series
\begin{eqnarray}
{\rm Ad}(\exp(s\,Y_i))Y_j
=Y_j-s\,[Y_i,Y_j]+\frac{s^2}{2}\,[Y_i,[Y_i,Y_j]]-\cdots,
\end{eqnarray}
where $s$ is a parameter and $i,j=1,\cdots,4$. The adjoint
representations of ${\goth g}_4$ is listed in Tables 2; it
consists the separate adjoint actions of each element of ${\goth
g}_4$ on all other elements.
\begin{table}
\centering{\caption{Commutator table}}\label{table:1}
\vspace{-0.3cm}\begin{eqnarray*}\hspace{-0.75cm}\begin{array}{l l
l l l} \hline
  [\,,\,]    &\hspace{3.1cm} Y_1      &\hspace{3.75cm} Y_2  &\hspace{3.75cm} Y_3      &\hspace{3.75cm} Y_4    \\ \hline
  Y_1        &\hspace{3.1cm} 0        &\hspace{3.75cm} 0    &\hspace{3.75cm} 2\,Y_1 &\hspace{3.75cm} 0    \\
  Y_2        &\hspace{3.1cm} 0        &\hspace{3.75cm} 0    &\hspace{3.75cm} Y_2    &\hspace{3.75cm} 0    \\
  Y_3        &\hspace{3.1cm} -2\,Y_1  &\hspace{3.75cm}-Y_2  &\hspace{3.75cm} 0      &\hspace{3.75cm} 0    \\
  Y_4        &\hspace{3.1cm} 0        &\hspace{3.75cm} 0    &\hspace{3.75cm} 0      &\hspace{3.75cm} 0    \\
  \hline
 \end{array}\end{eqnarray*}
\end{table}
\begin{table}
\centering{\caption{Adjoint table}}\label{table:2}
\vspace{-0.35cm}\begin{eqnarray*}\hspace{-0.75cm}\begin{array}{l l
l l l} \hline
  {\rm Ad}  &\hspace{3cm} Y_1      &\hspace{3.5cm} Y_2  &\hspace{3.25cm} Y_3           &\hspace{3.5cm} Y_4    \\ \hline
  Y_1       &\hspace{3cm} Y_1      &\hspace{3.5cm} Y_2  &\hspace{3.25cm} Y_3-2\,s\,Y_1 &\hspace{3.5cm} Y_4    \\
  Y_2       &\hspace{3cm} Y_1      &\hspace{3.5cm} Y_2  &\hspace{3.25cm} Y_3-s\,Y_2    &\hspace{3.5cm} Y_4    \\
  Y_3     &\hspace{3cm}e^{2s}\,Y_1 &\hspace{3.5cm}e^{s}\,Y_2 &\hspace{3.25cm} Y_3      &\hspace{3.5cm} Y_4    \\
  Y_4       &\hspace{3cm} Y_1      &\hspace{3.5cm} Y_2
  &\hspace{3.25cm} Y_3 &\hspace{3.5cm} Y_4    \\
  \hline
 \end{array}\end{eqnarray*}
\end{table}
{\thm An optimal system of one-dimensional Lie subalgebras of
nonlinear fin equation (\ref{eq:1}) is provided by those generated
by
\begin{eqnarray}\hspace{-0.7cm}\begin{array}{rlrl}
1)&A^1=Y_1=\partial_t, \hspace{2cm}  &4)&A^4=Y_4=e^{-u}\,E\,\partial_E,  \\
2)&A^2=Y_2=\partial_x, \hspace{2cm}  &5)& A^5= \alpha\,Y_1 + Y_2 = \alpha\,\partial_t+\partial_x, \\
3)& A^3=Y_3=2\,t\,\partial_t + x\,\partial_x - 2\,h\,\partial_h,
\hspace{2cm} &6) & A^6 = \beta\,Y_3 + Y_4=2\,\beta\,t\,\partial_t
+ \beta\,x\,\partial_x-2\,\beta\,h\,\partial_h +
e^{-u}\,E\,\partial_E, \label{eq:16}
\end{array}\end{eqnarray}
for nonzero constants $\alpha,\beta$.}\\

\noindent{\bf Proof.} Let ${\goth g}_4$ is the symmetry algebra of
Eq.~(\ref{eq:1}) with adjoint representation determined in Table 2
and
\begin{eqnarray}
Y=a_1\,Y_1+a_2\,Y_2+a_3\,Y_3+a_4\,Y_4
\end{eqnarray}
is a nonzero vector field of ${\goth g}$. We will simplify as many
of the coefficients $a_i$ as possible through proper adjoint
applications on $Y$. We follow our aim in the below easy cases.
\begin{description}
\item[Case 1]
At first, assume that $a_4\neq 0$. Scaling $Y$ if necessary, we
can consider $a_4$ to be 1 and so follow the problem with
\begin{eqnarray}
Y=a_1\,Y_1+a_2\,Y_2+a_3\,Y_3+Y_4.
\end{eqnarray}
According to Table 2, if we act on $Y$ by ${\rm
Ad}(\exp(\frac{1}{2}\,a_1\,Y_1))$, the coefficient of $Y_1$ can be
vanished:
\begin{eqnarray}
Y'=a_2\,Y_2+a_3\,Y_3+Y_4.
\end{eqnarray}
Then we apply ${\rm Ad}(\exp(a_2\,Y_2))$ on $Y'$ to cancel the
coefficient of $Y_2$:
\begin{eqnarray}
Y''=a_3\,Y_3+Y_4.
\end{eqnarray}
\item[Case 1a]
If $a_3\neq 0$ then we can not simplify the coefficient of $Y_3$
to be either $+1$ or $-1$. Thus any one--dimensional subalgebra
generated by $Y$ with $a_3, a_4\neq 0$ is equivalent to one
generated by $\beta\,Y_3 + Y_4$ which introduce parts {\it 6)} of
the theorem for constant $\beta\neq 0$.
\item[Case 1b]
For $a_3=0$ we can see that each one--dimensional subalgebra
generated by $Y$ is equivalent to $\langle\,Y_4\,\rangle$.
\item[Case 2]
The remaining one--dimensional subalgebras are spanned by vector
fields of the form $Y$ with $a_4=0$.
\item[Case 2a]
If $a_3\neq 0$ then by scaling $Y$, we can assume that $a_3=1$.
Now by the action of ${\rm Ad}(\exp(\frac{1}{2}\,a_1\,Y_1))$ on
$Y$, we can cancel the coefficient of $Y_1$:
\begin{eqnarray}
\bar{Y}=a_2\,Y_2+Y_3.
\end{eqnarray}
Then by applying ${\rm Ad}(\exp(a_2\,Y_2))$ on $\bar{Y}$ the
coefficient of $Y_2$ can be vanished and we tend to part {\it 3)}
of the theorem.
\item[Case 2b]
Let $a_3=0$ then  $Y$ is in the form
\begin{eqnarray}
\widehat{Y}=a_1\,Y_1+a_2\,Y_2.
\end{eqnarray}
Suppose that $a_2\neq 0$ then if necessary we can let it equal to
$1$. the simplest possible form of $Y$ is equal to
$\widehat{Y}=a_1\,Y_1+Y_2$ after taking $a_2=1$.
\item[Case 2b-1]
Let $a_1$ be nonzero. In this case we can not make the coefficient
of $Y_1$ in $\widehat{Y}$ more simpler and find {\it 5)} section
of the theorem.
\item[Case 2b-2]
If $a_1$ is zero then by scaling we can make the coefficient of
$\widehat{Y}$ equal to $1$. Hence this case suggests part {\it
2)}.
\item[Case 2c]
Finally if in the latter case $a_2$ be zero, then no further
simplification is possible and then $Y$ is reduced to ${\it 1)}$.
\end{description}

There is not any more possible case for studying and the proof is
complete.\hfill\ $\Box$
\\[3mm] \mbox{ }\hspace{3mm}    The coefficients $E, h$ of
Eq.~(\ref{eq:1}) resp. depend on the variables $u,x$. Therefore,
we take their optimal system's projections on the space $(x, u, E
, h)$. The nonzero in $x-$axis or $u-$axis projections of
(\ref{eq:16}) are
\begin{eqnarray}\begin{array}{rlrl}
1)&Z^1=A^2=A^5=\partial_x, \hspace{1cm} &3)&Z^3=A^4=e^{-u}\,E\,\partial_E,  \\
2)&Z^2=A^3=x\,\partial_x - 2\,h\,\partial_h, \hspace{1cm} &4)& Z^4
=A^6 = \beta\,x\,\partial_x-2\,\beta\,h\,\partial_h +
e^{-u}\,E\,\partial_E.\label{eq:16-1}
\end{array}\end{eqnarray}
From paper 7 of \cite{Ib3} we conclude that
{\prop{\label{a}} Let ${\goth g}_m:=\langle\,Y_i: i = 1, \cdots, m
\,\rangle$ be an $m$--dimensional algebra. Denote by $A^i\, (i =
1, \cdots, s,\, 0<s\leq m,\, s \in {\Bbb N})$ an optimal system of
one--dimensional subalgebras of ${\goth g}_m$ and by $Z^i\, (i =
1, \cdots, t,\, 0<t\leq s,\, t\in {\Bbb N})$ the projections of
$A^i$, i.e., $Z^i = {\rm pr}(A^i)$. If equations
\begin{eqnarray}
F = F(u),\hspace{0.75cm} g =g(x),
\end{eqnarray}
are invariant with respect to the optimal system $Z^i$ then the
equation
\begin{eqnarray}
u_t = (F(u)\,u_x)_x +  g(x)\,u,\label{eq:17}
\end{eqnarray}
admits the operators $X^i=$ projection of $A^i$ on $(t, x, u)$.}\\
{\prop{\label{b}} Let Eq.~(\ref{eq:17}) and the equation
\begin{eqnarray}
u_t = (\bar{F}(u)\,u_x)_x +  \bar{g}(x)\,u,\label{eq:18}
\end{eqnarray}
be constructed according to Proposition \ref{a} via optimal
systems $Z^i$ and $\bar{Z}^i$ resp. If the subalgebras spanned on
the optimal systems $Z^i$ and $\bar{Z}^i$ resp. are similar in
${\goth g}_m$, then Eqs.~(\ref{eq:17}) and (\ref{eq:18}) are
equivalent with respect to the equivalence group $G_m$ generated
by ${\goth g}_m$.}
\\[3mm] \mbox{ }\hspace{3mm}    Now by applying Propositions \ref{a}
and \ref{b} for the optimal system (\ref{eq:16-1}), we want to
find all nonequivalent equations in the form of Eq.~(\ref{eq:1})
admitting ${\mathcal E}$--extensions of the principal Lie algebra
${\goth g}_{{\mathcal E}}$, by one dimension, i.e, equations of
the form (\ref{eq:1}) such that they admit, together with the one
basic operators (\ref{eq:1-1}) of ${\goth g}_1$, also a second
operator $X^{(2)}$. In each case which this extension occurs, we
indicate the corresponding coefficients $E , h$ and the additional
operator $X^{(2)}$.
\\[3mm] \mbox{ }\hspace{3mm}    We perform the algorithm passing
from operators $Z^i\,(i=1,\cdots,4)$ to $E, h$ and $X^{(2)}$ via
the following example.
\\[3mm] \mbox{ }\hspace{3mm}   Let consider the vector field
\begin{eqnarray}
Z^4=-x\,\frac{\partial}{\partial x}+
2\,h\,\frac{\partial}{\partial h}+
e^{-u}\,E\,\frac{\partial}{\partial E},
\end{eqnarray}
then the characteristic equation corresponding to $Z^4$ is
\begin{eqnarray}
\frac{dx}{\beta\,x}=\frac{dh}{2\beta\,h}=\frac{dE}{e^{-u}\,E},
\end{eqnarray}
which determines invariants. Invariants can be taken in the
following form
\begin{eqnarray}
I_1=u, \hspace{0.5cm} I_2= h^{1/2}\,x, \hspace{0.5cm}
I_3=x^{-e^{-u}/\beta}\,E. \label{eq:19}
\end{eqnarray}

In this case there are no invariant equations because the
necessary condition for existence of invariant solutions (see
\cite{Ov}, Section 19.3) is not satisfied, i.e., invariants
(\ref{eq:19}) cannot be solved with respect to $E$ and $h$ since
$I_3$ is not an invariant function of $I_1$ and $I_2$ to derive a
function in term of $u$ for $E$.
\\[3mm] \mbox{ }\hspace{3mm}    Considering $Z^2$ we have the below
characteristic equation
\begin{eqnarray}
\frac{dx}{-x}=\frac{dh}{2h},
\end{eqnarray}
This equation suggest the following invariants
\begin{eqnarray}
I_1=u, \hspace{0.5cm} I_2= h^{1/2}\,x, \hspace{0.5cm}
I_3=E.\label{eq:20}
\end{eqnarray}
From the invariance equations we can write
\begin{eqnarray}
I_2=\phi(I_1), \hspace{1cm} I_3= \psi(I_1),
\end{eqnarray}
provided that $\phi(I_1)$ is a function in term of variable $u$
and $\psi(I_1)$ an arbitrary function. The first condition occurs
when $\phi(I_1)=c$ for constant $c$. It results in the forms
\begin{eqnarray}
E=\phi(u),\hspace{1cm} h=\Big(\frac{c}{x}\Big)^2.
\end{eqnarray}
From Proposition \ref{a} applied to the operator $Z^4$ we obtain
the additional operator
\begin{eqnarray}
X^{(2)}= 2\,\beta\,t\,\frac{\partial}{\partial t} +
\beta\,x\,\frac{\partial}{\partial x}.
\end{eqnarray}
\\[3mm] \mbox{ }\hspace{3mm}   One can perform the algorithm for
other $Z^i$~s of (\ref{eq:16-1}) similarly. The preliminary group
classification of nonlinear fin equation (\ref{eq:1}) admitting an
extension ${\goth g}_2$ of the principal Lie algebra ${\goth g}_1$
is listed in Table 3.
\begin{table}
\centering{\caption{The result of the
classification}}\label{table:3} \vspace{-0.35cm}\begin{eqnarray*}
\hspace{-0.75cm}\begin{array}{l l l l l l} \hline
  N       &\hspace{2cm} Z     &\hspace{2.1cm} \mbox{Invariant}  &\hspace{2cm} \mbox{Equation}
  &\hspace{2cm} \mbox{Additional operator}\,X^{(2)} \\ \hline
  1       &\hspace{2cm} Z^1    &\hspace{2.1cm} u         &\hspace{2cm} u_t=(\phi(u)\,u_x)_x+c\,u
  &\hspace{2cm} \frac{\partial}{\partial x},\,\alpha\,\frac{\partial}{\partial t}+\frac{\partial}{\partial x} \\
  2       &\hspace{2cm} Z^2    &\hspace{2.1cm} u         &\hspace{2cm} u_t=(\phi(u)\,u_x)_x+(\frac{c}{x})^2\,u
  &\hspace{2cm} 2\,t\,\frac{\partial}{\partial t} + x\,\frac{\partial}{\partial x} \\
  \hline
\end{array}\end{eqnarray*}
\end{table}
\section{Conclusion}
Projective analysis as a new symmetry property of equations
$u_t=(E(u)\,u_x)_x + h(x)\,u$ rather than previous results on this
equation \cite{PS,VP}, is carried out exhaustively. Also,
equivalence classification is given of the equation admitting an
extension by one of the principal Lie algebra of the equation. The
paper is one of few applications of a new algebraic approach to
the problem of group classification: the method of preliminary
group
classification. Derived results are summarized in Table 3. %

\end{document}